\input epsf

\magnification=1200


\hsize=125mm             
\vsize=195mm             
\parskip=0pt plus 1pt    
\clubpenalty=10000       
\widowpenalty=10000      
\frenchspacing           
\parindent=8mm           

\let\txtf=\textfont
\let\scrf=\scriptfont
\let\sscf=\scriptscriptfont
\font\frtnrm =cmr12 at 14pt
\font\tenrm  =cmr10
\font\ninerm =cmr9
\font\sevenrm=cmr7
\font\fiverm =cmr5

\txtf0=\tenrm
\scrf0=\sevenrm
\sscf0=\fiverm

\def\rm{\fam0 \tenrm}

\font\frtnmi =cmmi12 at 14pt
\font\tenmi  =cmmi10
\font\ninemi =cmmi9
\font\sevenmi=cmmi7
\font\fivemi =cmmi5

\txtf1=\tenmi
\scrf1=\sevenmi
\sscf1=\fivemi

 \def\oldstyle{\fam1 \tenmi}

\font\tensy  =cmsy10
\font\ninesy =cmsy9
\font\sevensy=cmsy7
\font\fivesy =cmsy5

\txtf2=\tensy
\scrf2=\sevensy
\sscf2=\fivesy

\font\tenit  =cmti10
\font\nineit =cmti9
\font\sevenit=cmti7
\font\fiveit =cmti7 at 5pt

\txtf\itfam=\tenit
\scrf\itfam=\sevenit
\sscf\itfam=\fiveit

\def\it{\fam\itfam\tenit}
\font\tenbf  =cmb10
\font\ninebf =cmb10 at  9pt
\font\sevenbf=cmb10 at  7pt
\font\fivebf =cmb10 at  5pt

\txtf\bffam=\tenbf
\scrf\bffam=\sevenbf
\sscf\bffam=\fivebf

\def\bf{\fam\bffam\tenbf}

\newfam\msbfam       
\font\tenmsb  =msbm10
\font\sevenmsb=msbm7
\font\fivemsb =msbm5

\txtf\msbfam=\tenmsb
\scrf\msbfam=\sevenmsb
\sscf\msbfam=\fivemsb

\def\msb{\fam\msbfam\tenmsb}
\def\Bbb#1{{\msb #1}}

\newfam\scfam
\font\tensc  =cmcsc10

\txtf\scfam=\tensc

\def\sc{\fam\scfam\tensc}


\def\frtnmath{%
\txtf0=\frtnrm        
\txtf1=\frtnmi         
}

\def\frtnpoint{%
\baselineskip=16.8pt plus.5pt minus.5pt%
\def\rm{\fam0 \frtnrm}%
\def\oldstyle{\fam1 \frtnmi}%
\everymath{\frtnmath}%
\everyhbox{\frtnrm}%
\frtnrm }


\def\ninemath{%
\txtf0=\ninerm        
\txtf1=\ninemi        
\txtf2=\ninesy        
\txtf\itfam=\nineit      
\txtf\bffam=\ninebf      
}

\def\ninepoint{%
\baselineskip=10.8pt plus.1pt minus.1pt%
\def\rm{\fam0 \ninerm}%
\def\oldstyle{\fam1 \ninemi}%
\def\it{\fam\itfam\nineit}%
\def\bf{\fam\bffam\ninebf}%
\everymath{\ninemath}%
\everyhbox{\ninerm}%
\ninerm }



\def\text#1{\hbox{\rm #1}}

\def\cite#1{{\uppercase{#1}}}
\def\ref#1{{\uppercase{#1}}}
\def\label#1{{\uppercase{#1}}}
\def\br{\hfill\break} 



\def\topmatter{\null\firstpagetrue\vskip\bigskipamount}
\def\endtopmatter{\vskip2\bigskipamount}

\def\title#1{%
\vbox{\raggedright\frtnpoint
\noindent #1\par}
\vskip 2\bigskipamount}        


\def\shorttitle#1{\rightheadtext={#1}}              

\newif\ifThanks
\global\Thanksfalse

\def\author#1{\begingroup\raggedright
\noindent{\sc #1\ifThanks$^*$\else\fi}\endgroup
\leftheadtext={#1}\vskip \bigskipamount}

\def\endabstract{\endgroup}

\long\def\abstract#1\endabstract{\par
\begingroup\ninepoint\narrower
\noindent{\sc Abstract.\enspace}#1%
\vskip\bigskipamount\endabstract}

\def\section#1#2{\bigbreak\bigskip\begingroup\raggedright
\noindent{\bf #1.\quad #2}\nobreak
\medskip\endgroup\noindent\ignorespaces}

\def\proclaim#1{\medbreak\noindent{\sc #1.\enspace}\begingroup
\it\ignorespaces}
\def\endproclaim{\endgroup\bigbreak}

\def\remark#1{\medbreak\noindent{\sc Remark \enspace}
\begingroup\ignorespaces}
\def\endremark{\endgroup\bigbreak}


\def\qed{$\mathord{\vbox{\hrule\hbox{\vrule
\hskip5pt\vrule height5pt\vrule}\hrule}}$}

\def\demo#1{\medbreak\noindent{{#1}.\enspace}\ignorespaces}
\def\enddemo{\penalty-100\null\hfill\qed\bigbreak}

\newdimen\EZ

\EZ=.5\parindent

\newbox\itembox

\newdimen\ITEM
\newdimen\ITEMORG
\newdimen\ITEMX
\newdimen\BUEXE

\def\iteml#1#2#3{\par\ITEM=#2\EZ\ITEMX=#1\EZ\BUEXE=\ITEM
\advance\BUEXE by-\ITEMX\hangindent\ITEM
\noindent\leavevmode\hskip\ITEM\llap{\hbox
to\BUEXE{#3\hfil$\,$}}%
\ignorespaces}


\newif\iffirstpage\newtoks\righthead
\newtoks\lefthead
\newtoks\rightheadtext
\newtoks\leftheadtext
\righthead={\ninepoint\rm\hfill{\the\rightheadtext}\hfill\llap{\folio}}
\lefthead={\ninepoint\rm\rlap{\folio}\hfill{\the\leftheadtext}\hfill}
\headline={\iffirstpage\hfill\else
\ifodd\pageno\the\righthead\else\the\lefthead\fi\fi}
\footline={\iffirstpage\hfill\global\firstpagefalse\else\hfill\fi}

\leftheadtext={}
\rightheadtext={}


\def\Refs{\bigbreak\bigskip\noindent{\bf References}\medskip
\begingroup\ninepoint\parindent=40pt}
\def\endRefs{\par\endgroup}
\def\endref{}

\def\ref{\par}
\def\key#1{\item{\hbox to 30pt{[#1]\hfill}}}


\def\cline#1{\leftline{\hfill#1\hfill}}

\def\bR{\Bbb R}
\def\bC{\Bbb C}
\def\bN{\Bbb N}
\def\bZ{\Bbb Z}

\def\supp{\mathop{\hbox{\rm supp}}\nolimits}
\def\kernel{\mathop{\hbox{\rm kernel}}\nolimits}

\input epsf

\topmatter
\title{Generic immersions of curves, knots,\br
monodromy  and gordian number}
\author{Norbert A'Campo}
\shorttitle{Immersions and knots.}
\endtopmatter

\hfill \hbox{\it Dedicated to Rob Kirby on his 60th birthday}

\par
\noindent
{\bf Table of contents}\par\noindent
1. Introduction \br
2. The fibration of the link of a divide \br
3. The monodromy diffeomorphism \br
4. Examples, symplectic and contact properties \br
5. The gordian number of the link of a divide \br

\S 1. {\bf Introduction}\par
A divide $P$ is a generic relative immersion of a finite number of 
copies of 
the unit interval $(I,\partial I)$ in the unit disk $(D,\partial D).$ 
The image of each copy of the unit interval is called a branch of the divide
$P$. 
The link $L(P)$ of a divide $P$ is
$$L(P):=\{(x,u)  \in T(P) \mid  \|(x,u)\|=1 \} \subset S(T(\bR^2))=S^3,$$
where we used the following notations:
For a tangent vector $(x,u) \in T(\bR^2)(=\bR^4)$  of $\bR^2$ the point 
$x \in \bR^2$ represents its foot and the vector $u \in T_x(\bR^2)$ its linear
part. The unit 
sphere 
$S(T(\bR^2)):=\{(x,u) \in T(\bR^2) \mid \|(x,u)\|:=x_1^2+x_2^2+u_1^2+u_2^2=1\}$ 
should not be confused with the 
tangent circle bundle of $\bR^2$ and is homeomorphic to the $3$-sphere 
$S^3.$
Finally, $T(P) \subset T(D) \subset T(\bR^2)$ is the space of tangent 
vectors of the divide $P,$ where at a crossing point $s$ by definition 
the space 
$T_s(P)$ is the union of the two 
$1$-dimensional subspaces of $T_s(D),$ 
which are the tangent spaces of the local branches of $P$ passing through 
$s.$
 The link $L(P)$ is an embedding of a union of $r$ circles in $S^3,$ 
where $r$ is the number of branches of the divide $P.$ So, for  a 
divide $P$ consisting of one branch 
the link 
$L(P)$ is a classical knot. \br

A divide is called connected if the image of the immersion is a connected 
subset of the disk. The following  is the main theorem of this paper. 

\proclaim{Theorem}
The link $L(P)$ of a connected divide $P$ is a fibered link.
\endproclaim
 
The monodromy of the fibered link $L(P)$ of a connected divide is given 
by Theorem $2$ 
of Section $3$ in terms of the combinatorics of the 
underlying divide. Since it is very easy to give 
examples of connected divides we 
obtain a huge class of links, such that the complement admits a fibration 
over the circle 
and that the isotopy class of the monodromy diffeomorphism is explicitly 
known. The links of plane curve singularities belong to this class  
(see [AC1,AC2,AC3,G-Z]). In Section~$5$ we show that the 
gordian number of the link of a divide equals the number of 
crossing points of the divide. The figure eight knot does not 
belong to this class. Many knots of this class are hyperbolic, as we will 
see in a forthcoming paper. Theorem~$2$ is  
used in the proof of the main theorem of [AC3]. 

I would like to thank Michel Boileau, 
Yasha Eliashberg, Rainer Kaenders, Dieter Kotschick, Tom Mrowka
and Bernard Perron for helpful discussions during the preparation of this
paper.

\S 2. {\bf The fibration of the link of a divide}\par
A regular isotopy of a divide in the space of generic immersions does not 
change the isotopy type of its link. So, without loss of generality, we may 
choose a divide to be linear and orthogonal near 
its crossing points. 
For a connected 
divide $P \subset D,$ let $f_P:D \to \bR$ be a generic 
$C^{\infty}$ function, 
such that $P$ is its $0$-level and that each region has exactly one 
non-degenerate maximum or minimum and that each region, which 
meets the boundary, 
has  exactly one non-degenerate maximum or minimum on the 
intersection of the region with $\partial D.$ Such a function exists for a 
connected divide 
and is well defined up to sign and isotopy. In particular, there are no 
critical points of saddle type other than the crossing points of the 
divide. Moreover without
loss of generality, we may assume  that 
the function $f_P$ is quadratic and 
euclidean in a neighborhood of those of its critical 
points, that lie in the interior
of $D$, i.e. for euclidean 
coordinates $(X,Y)$ with center at a  critical point~$c$ of $f_P$, in a
neighborhood of $c$ we have 
the expression $f_P(X,Y)=f_P(c)+XY$, if $c$ 
is a saddle point, $f_P(X,Y)=f_P(c)-X^2-Y^2$ if $c$ is a local maximum, or 
$f_P(X,Y)=f_P(c)+X^2+Y^2,$ if $c$ is a local minimum. Further, we may also 
assume, 
that the function $f_P$ is linear in a 
neighborhood of every relative 
critical point on $\partial D,$ i.e. at a critical point $c \in \partial D$ 
of the function $f_P$ we have the expression 
$f_P(c+h)=f_P(c)+<h,c>$ or 
$f_P(c+h)=f_P(c)-<h,c>$, where we denote by $<\,,\,>$ the scalar product of
$\bR^2$. Let $\chi:D \to [0,1]$ be a 
positive $C^{\infty}$ function which equals zero outside of 
small neighborhoods where $f_P$ is quadratic and equals $1$ in some smaller 
neighborhood $U$ of the critical points of $f_P.$ Moreover, we choose the 
function $\chi$ to be rotational symmetric around 
each critical point, i.e. we assume that locally near each critical point 
the function $\chi$ depends 
only on the distance to the critical point. For $\eta \in \bR, \eta > 0$ let 
$\theta_{P,\eta}:S^3 \to \bC$ be given 
by: 

$$
\theta_{P,\eta}(x,u):=f_P(x)+i \eta\, df_P(x)(u)-
{1\over 2}\eta^2\chi(x)H_{f_P}(x)(u,u)
$$ 

\noindent Observe that the Hessian $H_{f_P}$ is locally 
constant in the neighborhood of the  critical points of $f_P$ where $f_P$ is
euclidean.
Let 
$\pi_{P,\eta}:S^3 \setminus L(P) \to S^1$ be defined 
by:

$$
\pi_{P,\eta}(x,u):=\theta_{P,\eta}(x,u)/ |\theta_{P,\eta}(x,u)|.
$$

\proclaim{Theorem $1$}
Let $P$ be a connected divide. For $\eta >0$ and 
sufficiently small, the map $\pi_P:=\pi_{P,\eta}$ is a fibration 
of the complement of $L(P)$ over $S^1.$
\endproclaim

\demo{ \bf Proof} There exists a regular product 
tubular neighborhood $N$ of $L(P),$ such that the map $\pi_{P,\eta}$ 
for 
any $1 \geq \eta > 0$ is on $N \setminus L(P)$  a
fibration over $S^1$, 
for which  near $L(P)$ the fibers look like the pages of a book 
near its back. It is crucial to 
observe that in the intersection of the link $L(P)$ with 
the support of the function 

$$(x,u)\in S^3 \mapsto \chi(x)\in \bR$$ 

\noindent the kernel of the Hessian of 
$\theta_{P,\eta}$ and the kernel of the differential of the map 

$$(x,u)\in S^3 \mapsto f_P(x)\in  \bR$$ 

\noindent coincide. For any $\eta > 0$, 
the map $\pi_{P,\eta}$ is regular at each point of
$U':=\{(x,u) \in S^3 \mid x \in U\}.$ There exists $\eta_0>0$ such that 
for 
any $\eta, 0<\eta<\eta_0,$ the map $\pi_{P,\eta}$ is regular on 
$S^3\setminus (N\cup U').$ Hence, due to the quadratic 
scaling,  for $\eta$ sufficiently small the map $\pi_{P,\eta}$  is  a 
submersion, so since already a fibration near 
$L(P),$ it is a fibration by a theorem of Ehresmann. 
\enddemo

\S 3. {\bf The monodromy diffeomorphism}\par
Let $P$ be a connected divide and let $\pi_P:S^3\setminus L(P) \to S^1$ 
be its fibration of Theorem $1.$ We will show how to read off geometrically
the fibers $\pi_P^{-1}(\pm 1)$. Two diffeomorphisms $S''_i,S''_{-i}$ 
between the 
fibers $\pi_P^{-1}(\pm 1)$
modified by half Dehn twists will after a suitable composition give the
monodromy.
For our construction we orient the disk $D$ by one of the 
possible orientations, which we 
think of as an orthogonal  complex structure $J:T(D) \to T(D).$ We 
start out with a description of the fiber $F_1:=\pi_P^{-1}(1)$ and 
at the same time 
of the fiber $F_{-1}:=\pi_P^{-1}(-1).$ Put

$$
P_+:=\{x \in D \setminus \partial D \mid f_P(x) > 0, df_P(x) \not=0 \}
$$ 

\noindent The level 
curves of $f_P$ define  a oriented foliation $F_+$ on $P_+$, where 
a tangent vector $u$ to a level of $f_P$ at $x \in P_+$ is oriented if 
$df_P(x)(Ju) > 0.$ Put 

$$
P_{+,+}:=\{(x,u) \in S^3 \mid x\in P_+,\  u \in T(F_+)\}
$$ 

\noindent and 

$$
P_{+,-}:=\{(x,u) \in S^3 \mid x\in P_+,\  u \in T(F_-)\},
$$ 

\noindent where $F_-$ is 
the 
foliation with the opposite orientation. 
Put 

$$
F_M:=\{(x,u) \in S^3 \mid x=M\}
$$  

\noindent for a maximum $M,$ and 

$$
F_m:=\{(x,u) \in S^3 \mid x=m \}
$$ 

\noindent for a minimum $m$ of 
$f_P.$  Put 

$$
F_{s,+}:=\{(x,u) \in S^3 \mid x=s,\  H_{f_P}(x)(u,u) < 0\}
$$

\noindent and 

$$
F_{s,-}:=\{(x,u) \in S^3 \mid x=s,\  H_{f_P}(x)(u,u) > 0\}
$$ 

\noindent for a crossing point 
$s$ of $P,$ which is also a saddle point of $f_P.$ 
Observe that the angle in between 
$u,v \in F_m$ or $u,v \in F_M$ is a natural distance function on $F_m$ or 
$F_M,$ which allows us to identify $F_m$ and $F_M$ with a circle. Finally, put 

$$
\partial D_+:=\{x \in \partial D \mid f_P(x) > 0\}
$$

Let $p_{\bR}:S^3 \to D$ be the projection $(x,u) \mapsto x.$
The projection $p_{\bR}$ maps each of the 
sets $P_{+,+}$ and $P_{+,-}$ 
homeomorphically to $P_+.$ The sets $F_m$ or $F_M$  are homeomorphic 
to $S^1,$ if $M$ or $m$ is a maximum or minimum of $f_P$ respectively, and 
the sets $F_{s,\pm}$ 
are homeomorphic to a disjoint union of two  open 
intervals if $s$ is a crossing point of $P.$ The set $\partial D_+$ is 
homeomorphic to a disjoint union of 
open intervals. 
We have that $F_1$ and $F_{-1}$ are disjoint unions of these sets:

$$
F_1= P_{+,+} \cup P_{+,-} \cup \partial D_+\cup \bigcup_{s \in  P} F_{s,+} 
\cup \bigcup_{M \in  P_+} F_M
$$

\noindent and accordingly, with the obvious changes of signs:

$$
F_{-1}= P_{-,+} \cup P_{-,-} \cup \partial D_-\cup \bigcup_{s \in P} F_{s,-} 
\cup \bigcup_{m \in P_-} F_m
$$

\noindent In fact, for $(x,u) \in P_{+,+}\cup P_{+,-}$ we 
have $\theta_P(x,u) \in 
\bR_{>0}$ since 

$$
\theta_P(x,u):=f_P(x)+i \eta\ df_P(x)(u)-{1\over 2}\eta^2\chi(x)H_{f_P}(x)(u,u),
$$

\noindent where $f_P(x)>0,\ df_P(x)(u)=0,\ \chi(x)H_{f_P}(u,u)\leq 0.$  Hence, 
$P_{+,+}\cup P_{+,-}$ is an open and dense subset in $F_1.$ Forming 
the closure of $P_{+,+}\cup P_{+,-}$ in $F_1$ leads to
the following combinatorial 
description of the above decomposition. First, we add to 
the open surface $F_1$ its boundary and get 

$$
\bar F_1:=F_1 \cup L(P)
$$

\noindent Let $R$ be a connected component of $P_+.$ The inverse image
$p_{\bR}^{-1}(R) \cap \bar F_1$ in $\bar F_1$ are two disjoint open cells or
cylinders 
$R_+ \subset P_{+,+}$ and $R_- \subset P_{+,-}$ 
which are in fact subsets of $F_1.$ The closure of 
$R_+$ in  $\bar F_1$ is a surface $\bar R_+$ with boundary  and corners. 
The set $F_M$ is a common boundary component without corners of 
$\bar R_+$ and $\bar R_-$ if $M$ is a maximum in $R.$ If there is 
no maximum in $R$ the closures  $\bar R_+$ and $\bar R_-$ meet along the 
component of $\partial D_+$ which lies in the 
closure of $R.$ Let $S,R$ be  connected components of $P_+$ such that the 
closures of $R$ and $S$ have a crossing point 
$s$ in common. The closures of $R_+$ and $S_-$ in $F_1$ meet along one of 
the components of $F_{s,+}$ and the closures of 
$R_-$ and $S_+$ in $F_1$ meet along the other component of $F_{s,+}.$ The 
closure of 
$F_{s,+} \cap \bar R_+$ in 
$\bar R_+$ intersects $L(P)$ in $2$ corners, that are  also corners of 
the closure of $F_{s,+} \cap \bar S_-$ in 
$\bar S_-$ (see Fig. 1). Notice that the foliation $F_+$ on $P_+$ does 
not lift to a
foliation, which extends to an oriented foliation on $F_1$.

\midinsert 
\cline{\epsffile{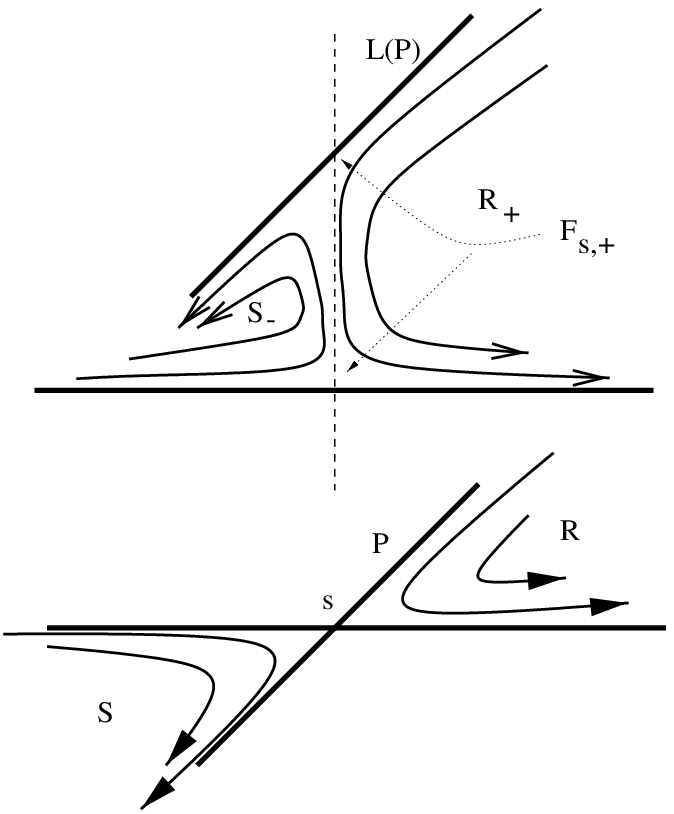}} 
\medskip 
\centerline{Figure 1: Gluing of the lifts of R 
with + and S with - foliation to $F_1$.}
\endinsert

Now we will work out the fibers $F_i:=\pi_P^{-1}(i)$ and 
$F_{-i}:=\pi_P^{-1}(-i).$ First observe that 
$F_i$ and $F_{-i}$ are projected  to a subset of 
$P \cup \supp(\chi)$ by  $p_{\bR}$. Put

$$
F_{i,P}:=\{(x,u) \in S^3 \mid x \in P,\  \chi(x)=0,\ df_P(x)(u)>0 \}.
$$

\noindent For a crossing point $c$ of $P$ we put

$$
F_{i,c}:=\{(x,u) \in S^3 \mid \chi(x) > 0,\  df_P(x)(u) > 0 ,\  
f_P(x)- {1\over 2}\eta^2\chi(x)H_{f_P}(c)(u,u)=0\}.
$$

\noindent In order to get nice sets it is necessary to choose 
a nice bump function 
$\chi.$
The set 
$F_{i,P}\cup F_{i,c}$ is an open and dense subset in $F_i$ and forming 
its closure in $\bar F_i:=F_i \cup L(P)$ leads to a 
combinatorial description of $F_i.$

Our next goal is the description of the monodromy diffeomorphism. 
We will use the integral 
curves of  the distribution $J(\kernel(df_P)),$ which 
pass through the crossing points of the divide $P.$ In a connected component 
$R$ of $D \setminus P$, those 
integral curves of  
$J(\kernel(df_P))$ meet at the critical point of $f_P$ in the component 
$R$ with distinct tangents 
or
go to distinct points of $\partial D.$ 

\midinsert 
\cline{\epsffile{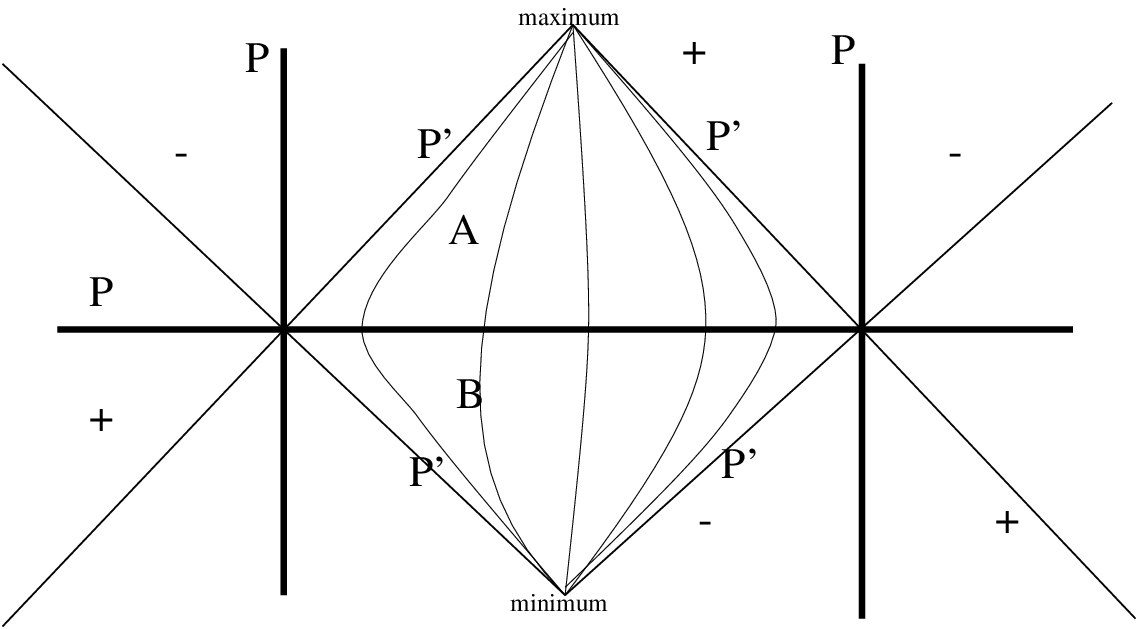}} 
\medskip 
\centerline{Figure 2: Two tiles with the $J(\kernel(df_P))$ foliation.}
\endinsert

We denote by $P'$ the union of the integral curves of $J(\kernel(df_P)),$ 
which pass through the crossing points of $P.$ The complement in $D$ of the 
union  $P' \cup P \cup \partial D$ is a disjoint 
union of tiles, which are homeomorphic to open squares or triangles. We call 
a pair $(A,B)$ of tiles opposite, if $A\not= B$ and 
the closures of $A$ and $B$ in $D$ have a segment of $P$ in common. For an 
opposite pair of tiles $(A,B)$ let $A|B$ be the 
interior in $D$ of the union of the closures of $A$ and $B$ in $D$. The set 
is foliated by the levels of $f_P$ and also by the 
integral lines of the distribution $J(\kernel(df_P)).$ Both foliations are 
non-singular and 
meet in a $J$-orthogonal way (see Fig. 2).

Let $R$ and $S$ be the components of $D\setminus (P\cup \partial D),$ which 
contain $A$ and $B.$ Put

$$
F_{1,A|B}:=\{(x,u) \in F_1 \mid x \in A|B \}
$$

\noindent and

$$
F_{-1,A|B}:=\{(x,u) \in F_{-1} \mid x \in A|B \}
$$

\noindent The sets $F_{1,A|B}$ and $F_{-1,A|B}$ each have 
two connected components:

$$
F_{1,A|B}=F_{1,+,A|B} \cup F_{1,-,A|B}
$$

\noindent where

$$
F_{1,+,A|B}:=\{(x,u) \in F_1 \mid x \in A|B,\  df_P(Ju) > 0 \}
$$

$$
F_{-1,+,A|B}:=\{(x,u) \in F_{-1} \mid x \in A|B,\  df_P(Ju) > 0 \}
$$

\noindent and

$$
F_{1,-,A|B}:=\{(x,u) \in F_1 \mid x \in A|B,\  df_P(Ju) < 0 \}
$$

$$
F_{-1,-,A|B}:=\{(x,u) \in F_{-1} \mid x \in A|B,\  df_P(Ju) < 0 \}
$$

The closures of $F_{1,\pm,A|B}$ in $\bar F_1$ and of $F_{-1,\pm,A|B}$ in 
$\bar F_{-1}$ are  polygons with $6$ edges: let 
$M,c,c'$ be the vertices of the triangle $A$; the six edges of the closure 
$H$ of $F_{1,+,A|B}$ in $\bar F_1$ are 
$\{(x,u) \in H \mid x=M\},\ $ 
$\{(x,u) \in H \mid x \in [c,M]\},\ $ 
$\{(x,u) \in H \mid x=c\},\ $ $\{(x,u) \in H \mid x \in [c,c']\},\ $ 
$\{(x,u) \in H \mid x=c'\},\ $ $\{(x,u) \in H \mid x \in [c',M]\}$
where $[M,c]$ and $[M,c']$ are segments included in $P'$ and $[c,c']$ is a 
segment in $P.$
 
We will define two diffeomorphisms:

$$
S_{i,A|B}: F_{1,A|B} \to  F_{-1,A|B}
$$

\noindent and

$$
S_{-i,A|B}: F_{1,A|B} \to  F_{-1,A|B}
$$

\noindent To do so it is convenient to choose the function $f_P:D \to \bR$ such 
that the maxima are of value $1$ and the minima of value $-1.$ Moreover, we
modify the function $f_P$ at the boundary $\partial{D}$ such that
along each of the
integral lines of the foliation given by the distribution $J(\kernel(df_P))$ 
the function $f_P$ takes  
all values in
an interval $[-m,m]$ with $1 \geq m > 0$. The latter modification of $f_P$ is
useful if the tile $A$ or $B$ meets $\partial{D}$. We also 
need the 
rotations $J_{\theta}:T(D) \to T(D)$ about the angle $\theta \in [-\pi,\pi].$ 
Remember $J=J_{\pi/2}.$ 
The map $S_i$ acts as follows: for $(x,u) \in F_1$ with $x\in A|B$ let 
$y\in A|B$ be the point in the opposite tile on the integral line of the 
distribution $J(\kernel(df_P))$ with $f_P(x)=-f_P(y);$  
now we move $x$ to $y$ along
the integral curve $\gamma_{x,y}(t), t \in [f_P(x),f_P(y)]$  which connects 
$x$ and $y$ with the parameterization  
$f_P(\gamma_{x,y}(t))=t$; the vector 
$u$ will be moved along the  path

$$
(\gamma_{x,y}(t),U_{x,y}(t)):=(  \gamma_{x,y}(t),
s(x,t)( J_{\theta(x,t)}(|f_P(x)| u_{x,y}(t)/2)+ u_{x,y}(t) )  ),
$$

\noindent where $(\gamma_{x,y}(t),u_{x,y}(t)) \in S^3$ is 
the continuous vector field 
along $\gamma_{x,y}(t)$ such that $u_{x,y}(t)$ 
stays in the kernel of 
$df_P$ and $u_{x,y}(f_P(x))=U_{x,y}(f_P(x))=u$, where the rotation angle $\theta(x,t)$ 
at time $t$ is given by $\theta(x,t):={(|f_P(x)|-t)\pi \over 2|f_P(x)|}$ 
and where the stretching factor $s(x,t) \geq 1$ is chosen such that
$(\gamma_{x,y}(t),U_{x,y}(t)) \in S^3$ holds;
define

$$
S_i((x,u)):=(y,u_{x,y}(f_P(y)))=(y,U_{x,y}(f_P(y)))
$$

\noindent The definition of $S_{-i}$ 
is analogous, but uses rotations in the sense of $-J.$ 
The names $S_i$ or $S_{-i}$ indicate that the flow lines 
$(\gamma_{x,y}(t),u(t))$ pass through the fiber 
$F_i$ or  $F_{-i}$ respectively. The flow lines defining $S_i$ or $S_{-i}$ are 
different. However, the maps $S_i$ and 
$S_{-i}$ are equal. The system of paths $(\gamma_{x,y}(t),U_{x,y}(t)) \in S^3$ 
is local
near the link $L(P)$, i.e. for every neighborhood $N$ in $S^3$ of a point
$(x',u') \in L(P)$ there exists a neighborhood $M$ of $(x',u')$ in $S^3$ 
such that
each path $(\gamma_{x,y}(t),U_{x,y}(t))$ with $(x,u) \in F_1 \cap M$ stays in
$N$. It will follow that the flow lines of the monodromy 
vector field are  meridians of the link $L(P)$ in its neighborhood.

The partially defined diffeomorphisms $S_i$ and 
$S_{-i}$

$$S_i,S_{-i}:\bigcup_{A|B} F_{1,A|B} \to \bigcup_{A|B} F_{-1,A|B}$$

\noindent are obtained by gluing 
the maps  $S_{i,A|B}: F_{1,A|B} \to  F_{-1,A|B}$ and 
$S_{-i,A|B}: F_{1,A|B} \to  F_{-1,A|B}$ for all
opposite pairs of tiles $(A,B)$ with $A\subset P_+$. The gluing poses 
no problem 
since those  unions are disjoint, but  the diffeomorphisms $S_i$ and $S_{-i}$
do not extend continuously to $F_1.$ We will see that the 
discontinuities, which are the obstruction for extending $S_i$ and 
$S_{-i},$ can be 
compensated by a composition of  right half Dehn twists.

\midinsert 
\cline{\epsffile{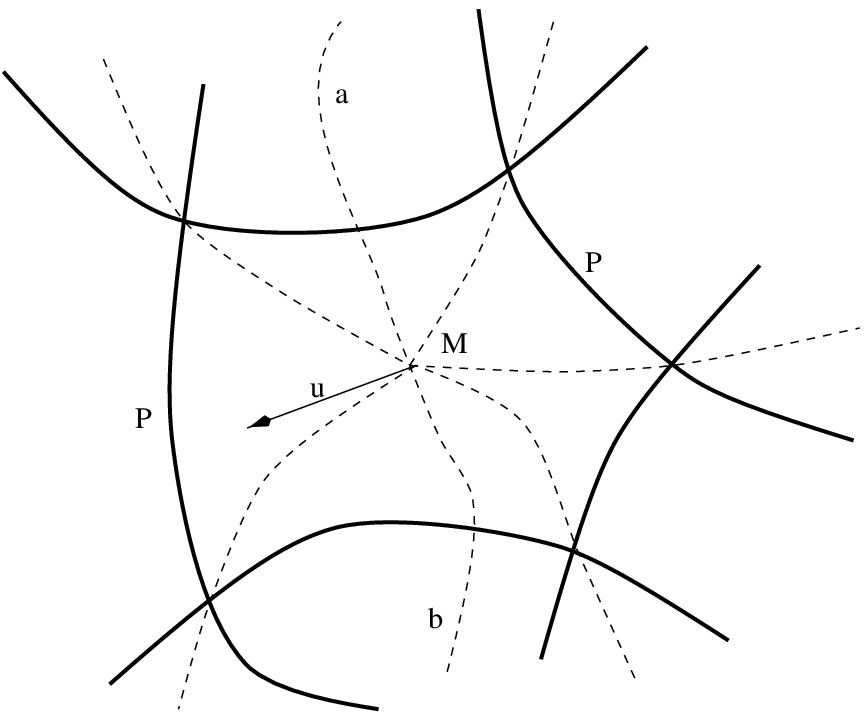}} 
\medskip 
\centerline{Figure 3: The discontinuity at $F_M.$}
\endinsert

At a maximum $M \in D$ of $f_P$ each vector $(M,u) \in S^3$ belongs to $F_1.$ 
Let $a$ and $b$ be the integral curves 
of $J(\kernel(df_P))$ with one endpoint at $M$ and  orthogonal to $u.$ We 
assume that neither $a$ nor $b$ passes through 
a crossing point of $P$ (see Fig. $3$) and that $a$ and $b$ belong to
different pairs of opposite tiles. A continuous extension of the 
maps $S_i$ or $S_{-i}$ has to map the vector $(M,u)$ 
to  two vectors based at the other 
endpoint of $a$ and $b.$ Since these endpoints 
differ in general, a continuous extension is impossible. 

In order to allow a continuous extension at the common endpoint of 
$a$ and $b$ we make a new surface $F'_1$ by 
cutting $F_1$ along the cycles $F_M,$ where $M$ runs through all the 
maxima of 
$f_P$ and by gluing back after a rotation of angle $\pi$ of each of the 
cycles  
$F_M.$ In the analogous manner, we make the surface $F'_{-1}$ in doing 
the half twist along $F_m,$ where $m$ runs through the minima of $f_P.$ The 
subsets $F_{1,A|B}$ do not meet the support of the 
half twists, so they are canonically again subsets of $F'_1,$ which we denote 
by $F'_{1,A|B}.$ Analogously, we have 
subsets $F'_{-1,A|B}$ in $F'_{-1}.$ 
A crucial observation is that the partially defined diffeomorphisms 
$$
S'_i,S'_{-i}:\bigcup_{A|B}F'_{1,A|B} \to \bigcup_{A|B}F'_{-1,A|B}
$$  
have less discontinuities, which are the obstruction for a continuous extension.
We denote by $a'$ and $b'$ the arcs on $F'_1,$ which correspond to 
the arcs $a$ and $b$ on $F_1.$ 
Indeed, the continuous extension  at the  end points of $a'$ and $b'$  
is now possible. 

Let $s$ be a crossing point of $P$ and let $I_{s,+}$ be the segment 
of $P',$ which passes through $s$ and lies in $P_+.$ The inverse 
image of $Z^{\circ}_s:=p_{\bR}^{-1}I_{c,+} \cap F_1$  is not a cycle, except 
if both endpoints of $I_{s,+}$ lie on $\partial D.$ If a maximum 
$M$ of $f_P$ is an endpoint of $I_{s,+},$ the inverse image 
$p_{\bR}^{-1}(M) \cap F_1$ consists of $2$ points on $F_M,$ 
which are antipodal. On the new surface $F'_1$ the inverse image 
$p_{\bR}^{-1}(I_{s,+}) \cap F'_{-1}$ is a cycle. 
An extension of 
$S'_i$ and $S'_{-i}$ will be discontinuous along this cycle 
(see Fig. $4$).
We now  observe that the partially defined diffeomorphisms $S'_i$ 
and 
$S'_{-i}$
have discontinuities along the cycle 
$p_{\bR}^{-1}(I_{s,+}) \cap F'_{-1}$, which can be compensated 
by half twists along the 
inverse images $p_{\bR}^{-1}(I_{s,-}) \cap F'_{-1},$ where $s$ 
runs through the crossing points of $P.$ 
Note that for a crossing point $s$ of $P$ the curve
$Z'_{s,-1}:=p_{\bR}{-1}(I_{s,-}) \cap F'_{-1}$ is in fact a simply 
closed curve on $F'_{-1}.$

\midinsert 
\cline{\epsffile{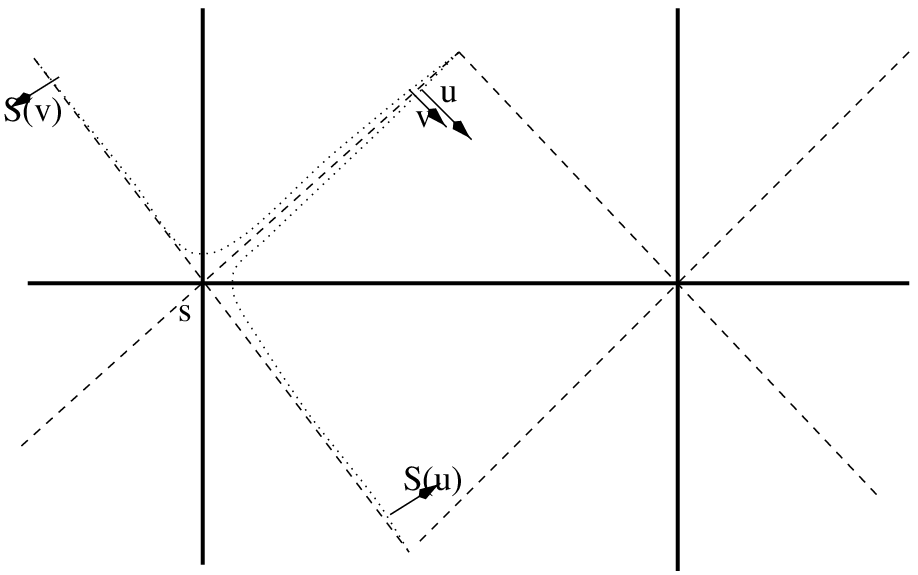}} 
\medskip 
\centerline{Figure 4: The discontinuity along $Z_s.$}
\endinsert

For a crossing point $s$ of the divide $P$ we now define a simply 
closed curve on $F_1,$ by putting:
$$
Z_s:=Z^{\circ}_s \cup \bigcup_{M \in \partial I_{s,+}} F_{s,M}
$$
where for an endpoint $M$ of $I_{s,+},$ which is a maximum of $f_P,$ 
the set 
$ F_{s,M}$ is the simple arc of $F_M,$ which connects the two points 
of $Z^{\circ}_s \cap F_M$ and contains an 
inward tangent vector of $I_{s,+}$ at $M.$ As we already have noticed 
the set $Z^{\circ}_s \cap F_M$ has only 
one element if $M \in \partial D,$ so we define $F_{s,M}:=\emptyset$ 
in that case.

We have the inclusion $F_m \subset F_{-1}.$ We now define the cycle 
$Z_m \subset F_1.$ Define for a minimum $m$ of $f_P$ the region
$$
B_m:= \bigcup_{A|B, m \in \bar{B}} F_{-1,A|B}
$$
Let $B_{m,\epsilon}$ be the level curve
$$
B_{m,\epsilon}:=\{(x,u) \in B_m \mid f_P(x)=-\epsilon\}
$$
For a small $\epsilon$ the set 
$$(S_i)^{-1}(B_{m,\epsilon} \cap \bigcup_{A|B} F_{-1,A|B})$$ 
is a union of copies of an open interval and is not a cycle but 
nearly a cycle. The unions closes up to a cycle by 
adding small segments which project to the integral lines through the 
crossing points of $P.$ We denote this cycle by 
$Z_m \subset F_1$.

We are now able to state the main theorem.

\proclaim{Theorem 2}
Let $P$ in $D$ be a connected divide. Let 
$\pi_P:S^3\setminus L(P) \to S^1$ be the fibration of Theorem $1.$ The counter 
clockwise monodromy of the fibration
$\pi_P$ is the composition of right Dehn twists 
$T:=T_- \circ T_. \circ T_+:F_1 \to F_1,$ where $T_-$ is 
the product of the right twists along $Z_m,$ $m$ running through the 
minima of $f_P,$ $T_.$ is the product of the right Dehn 
twists along  the cycles $Z_s,$ $s$ 
running through the crossing points of $P,$ and $T_+$ is the product 
of the right twists along $F_M,$  $M$ running through the maxima of $f_P.$
\endproclaim

Before giving the proof, we will define positive and negative half Dehn 
twists.
Let $X$ be an oriented surface and let $z$ be a simply closed curve on $X.$
Let $X'$ be the surface obtained from the surface $X$ by cutting $X$ along 
$z$ and by gluing back with a diffeomorphism of degree one. The surfaces 
$X$ and $X'$ are of course diffeomorphic. A minimal positive pair 
of Dehn twists from $X$ to $X'$ is a pair of 
diffeomorphisms $(p,q)$ from $X$ to $X'$ such that the following holds:

\item{(a)} The composition $q^{-1} \circ p:X \to X$ is a right Dehn twist 
with respect to the orientation of $X$ having 
the curve $z$ as core. In addition $p(z)=q(z)=z$ holds.

\item{(b)} There exists a regular bicollar neighbourhood $N$ of $z$ in $X$ 
such that  both $p$ and $q$ are the identity 
outside $N.$

\item{(c)} For some volume form $\omega$ on $N,$ which we think of as a 
symplectic structure,  we have 
$p^*\omega=q^*(\omega)=\omega$, and the sum of the Hofer  distances ([H-Z], see
Chap. 5,)
to the identity of the restrictions 
of $p$ and $q$ to  $N \setminus z$ is minimal. 

\noindent 
Minimal positive pairs of Dehn twists exist and are well defined up to isotopy. 
For a minimal positive pair $(p,q)$ of 
Dehn twists, the member $p$ is called positive or 
right and the member $q$ is called negative or left.

\demo{ \bf Proof of Theorem 2} We need to introduce one more surface. 
Let $F''_{-1}$ 
be the surface obtained from the surface $F'_{-1}$ by cutting $F'_{-1}$ 
along the cycles $Z'_{s,-1}$ and by gluing back after a half
twist along each $Z'_{s,-1},$ $s$ running through the 
crossing points of $P.$ We still have partially defined diffeomorphisms
$$
S''_i,S''_{-i}:\bigcup_{A|B}F'_{1,A|B} \to \bigcup_{A|B}F''_{-1,A|B}
$$
since the cutting was done in the complement of 
$\bigcup_{A|B}F'.$ By a direct inspection we see that the 
diffeomorphisms extend continuously to
$$
S''_i,S''_{-i}:F'_1 \to F''_{-1}
$$
Let 
$$(p_+,q_+):F_1 \to F'_1$$
$$(p_.,q_.):F''_{-1} \to F'_{-1}$$
$$(p_-,q_-):F'_{-1} \to F_{-1}$$
be minimal positive pairs of Dehn twists.
A direct inspection shows that the composition
$$
p_+ \circ S'_i \circ p_. \circ p_- \circ (q_+ \circ S'_{-i} \circ q_. \circ q_-)^{-1}:F_1 \to F_1
$$
is the monodromy of the fibration $\pi_P.$
This composition evaluates to
$$
T_- \circ T_. \circ T_+:F_1 \to F_1
$$
 
\enddemo
\goodbreak
\remark{\bf Remark $1$}
We list   special properties of the monodromy of links and knots of divides.
The number of Dehn twists of the above decomposition of the  monodromy 
equals the 
first betti number $\mu=2\delta - r +1$ of the fiber, and the total number of 
intersection points among the core curves of the involved Dehn twists 
is less then $5\delta.$ This means that the complexity of
the monodromy is bounded by a function of $\mu.$ For instance, the
coefficients of the Alexander polynomial of the link of a divide are 
bounded
by a quantity, which depends only on the degree of the 
Alexander polynomial. This observation suggests the
following definition for the complexity $C$ of an element of the mapping class
group $\phi$ of a surface: the minimum of 
the quantity $L+I$ over all decompositions 
as product of Dehn twists of $\phi,$ where $L$ is the number of factors and
$I$ is the number of mutual intersections of the core curves. We do not know
properties of this exhaustion of the mapping class group. Notice, that the
function $(\phi,\psi) \mapsto C(\psi^{-1}\circ \phi)\in \bN$ defines a left
invariant distance  on the mapping class group.

A crossing point of $P$ is 
``\`a quatre vents'', if the $4$ sectors of the complement 
of $P$ in $D$ meet the boundary 
of $D.$ If there are no crossing points \`a quatre vents, then none of 
the core curves of the twists involved in the decomposition of the monodromy 
does  separate the fiber, so the twists of the 
decomposition are all conjugated in the  orientation preserving mapping 
class group of the fiber.

It is easily seen that for any link of a divide the  
monodromy diffeomorphism and its inverse are conjugate by an
orientation reversing element in the mapping class group. In our previous
notations this conjugacy is given by the map $(x,u)\in F_1 \mapsto (x,-u) \in
F_1$, which moreover  realizes geometrically the symmetry of G. Torres [To]
$t^{\mu}\chi(1/t)=(-1)^{\mu}\chi(t)$
for the Alexander polynomial $\chi(t)$ of knots.

\endremark
  
\remark{\bf Remark $2$}
In fact the proof of Theorem $2$ shows that the fibration of the link of a 
connected divide $P$ can be filled with a singular fibration in the 
$4$-ball, which has $3$ singular fibers with only quadratic singularities, 
as in the case of a divide of the singularity of a complex plane curve. 
The filling has only two singular fibers if the 
function $f_P$ has no maxima 
or no minima. By this construction  
from a connected divide we obtain a 
contractible $4$-dimensional piece with a Lefschetz pencil. It is 
sometimes
possible 
to glue these local pieces and to get $4$-manifolds with a Lefschetz 
pencil.
    
\endremark

\S 4. {\bf Examples, symplectic and contact properties}\par

The figure eight is not the knot of a divide. The figure eight knot's 
complement fibers over the circle with as fiber the punctured torus and 
as monodromy the isotopy class of the 
linear diffeomorphism given by a matrix in $SL(2,\bZ)$ of trace $3.$ 
Such a matrix $M$ is not the product of two  
unipotent matrices, which are conjugate in $SL(2,\bZ)$ and the matrices $M$
and $M^{-1}$ are not conjugated by an integral matrix of determinant $-1$. 
So according 
to the remarks  of Section $3$,  the 
figure eight can not be the knot of a divide. A third argument to 
rule out the figure eight as the knot of a 
divide goes as follows. The first betti number of the fiber of the 
figure eight knot is $2.$  But only two 
connected divides
have  fibers with betti number $2$ and these two have monodromies with 
trace equal to $1$ or $2.$   

The connected sum of two divides $(D_1,P_1)$ and $(D_2,P_2)$ is done by 
making a boundary connected sum of $D_1$ and $D_2$ such that a boundary 
point of $P_1$ matches with a boundary point of $P_2.$ For 
divides 
with one branch we have the
formula:

$$
L(P_1\#P_2)=L(P_1)\#L(P_2)
$$

The Theorems $1$ and  $2$ remain true for generic immersions of 
disjoint unions of intervals and circles in the 
$2$-disk. It is also possible to start with a generic immersion of a
1-manifold $I$ in an oriented compact connected surface with boundary $S.$
The pair $(S,I)$ defines a link $L(S,I)$ in the 3-manifold 
$M_S:=T^+(S)/\hbox{\rm zip}$, where $T^+(S)$ is the
space of oriented tangent directions of the surface $S$ and 
where $\hbox{\rm zip}$ is the
identification relation, which identifies $(x,u),(y,v)\in T^+(S)$ 
if and only if 
$x=y\in \partial S$ or if $(x,u)=(y,v).$ In order to get a fibered 
link,  the topological pair  
$(R,R \cap \partial S)$ has to be contractible for each connected 
component $R$
of $S\setminus I$ and moreover, the complement $S\setminus I$ has to allow a
chess board coloring in positive and negative regions. 
The proofs do not need any modification.

A  relative immersion $i:I \to D$ of a copy of $[0,1]$ in $D,$ such 
that at selftangencies the
velocities are with opposite orientation, defines an embedded and 
oriented arc $I'$ in $S^3$ by putting:

$$
I':=\{(x,u) \in S^3 \mid x \in i(I), (di^{-1})_x(u) \geq 0\}
$$

Let $j:I \to D$ be a relative immersion with only transversal crossings 
and opposite selftangencies, such that
the endpoints of $i$ and $j$ are tangent with opposite orientations and that 
all tangencies of $i$ and $j$ are generic and have opposite orientations. 
The union $I' \cup J'$
is the oriented knot of the pair $(i,j).$ A divide $P$ 
defines pairs $(i_P,j_P)$ of relative immersions 
with opposite orientations by taking both orientations. Those pairs 
$(i_P,j_P)$ have a special $2$-fold symmetry. For instance 
the complex conjugation 
realizes this $2$-fold 
symmetry for a  divide, which arises as a real deformation of a real
plane curve singularity. It is in\-teresting to observe that this symmetry 
acts on $F_1$ with as
fixed point set the intersection $D\cap F_1$, which is a collection of $r$
disjointly embedded arcs in $F_1.$ The quotient of $\bar{F_1}$ by the symmetry
is an orbifold surface with exactly $2r$ boundary 
${\pi \over 2}$-singularities. 
Any link of singularity of a plane curves can be obtained as the link of a 
divide (see [AC3]). It is an 
interesting problem to characterize links of 
singularities among links of divides.

We finish this section with some symplectic and contact properties.
The link of a divide is transversal to the standard contact structure in the
$3$-sphere. This can be seen explicitly by the following computation, where
we use the multiplication of quaternions. Let $P$ be a divide in the unit
disk. We assume that the part of $P$, which lies in the collar of
$\partial{D}$  
with inner radius $1 \over \sqrt{2}$,
consists of radial line segments. We think of 
the branches of $P$ as parametrized
curves $\gamma(t)=(a(t),b(t)),-A \leq t \leq A,$ 
where the parameter speed is adjusted such that
$a^2+b^2+\dot{a}^2+\dot{b}^2=1.$ To the branch $\gamma$ correspond two arcs
$\Gamma^+$ and $\Gamma^-$ on the sphere of quaternions of unit length:

$$
\Gamma^+(t):=a(t)-\dot{a}(t)i+b(t)j+\dot{b}(t)k
$$

$$
\Gamma^-(t):=a(-t)+\dot{a}(-t)i+b(t)j-\dot{b}(-t)k
$$

The left invariant speed of $\Gamma^+$ at time $t$ is

$$
v(t):=\Gamma^+(t)^{-1}{d \over dt}\Gamma^+(t)
$$

We have

$$
v=a\dot{a}+\dot{a}\ddot{a}+b\dot{b}+\dot{b}\ddot{b}+
[-a\ddot{a}+\dot{a}^2-b\ddot{b}+\dot{b}^2]i+v_jj+v_kk
$$
The coefficient 
$v_0:=a\dot{a}+\dot{a}\ddot{a}+b\dot{b}+\dot{b}\ddot{b}$ vanishes, 
since ${d\over dt}\Gamma(t)$ is perpendicular to $\Gamma(t),$ and hence we 
can
rewrite the coefficient $v_i$ of $i$ in $v$ as

$$
v_i=-a\ddot{a}+\dot{a}^2-b\ddot{b}+\dot{b}^2=<(a+\dot{a},b+\dot{b})\mid
(\dot{a},\dot{b})>
$$

Outside of the collar neighborhood of $\partial{D}$ 
we have $v_i > 0,$ since $a^2+b^2 < 1/2 <
\dot{a}^2+\dot{b}^2.$ In the collar we also have $v_i > 0$ by a direct
computation. Since the left
invariant contact structure on the unit sphere in the skew field 
of the quaternions 
is given by the span of the tangent vectors 
$j$ and $k$ at the point $1$, we conclude that $\Gamma^{+}$ with its
orientation is in the positive sense transversal to the left invariant 
contact structure $S^3$. 

For the link of a divide we now will construct 
a polynomial, hence symplectic, spanning surface in the
$4$-ball. For $\lambda \in \bR, \lambda > 0$ put 

$$
B_{\lambda}:=\{p+ui \in \bC^2 \mid p,u \in \bR^2,||p||^2 + \lambda^{-2} ||u|| \leq 1\}
$$

We have $B_{\lambda} \cap \bR^2 = D$ and $B_{\lambda}$ is a strictly
holomorphically convex domain with smooth boundary in $\bC^2$. The map 
$(p,u) \mapsto ((p,u/{\lambda})$ identifies $\partial{B_{\lambda}}$ with
the unit $3$-sphere of $\bC^2$.

\proclaim{\bf Theorem 3}
Let $P$ be a connected divide in the disc $D$ with $\delta$ double points and
$r$ branches.
There exist $\lambda > 0, \eta  > 0$
and
there exists a polynomial function $F:B_{\lambda} \to \bC$ with the following
properties:

\item{(a)} the function $F$ is real, i.e. $F(\overline{p+ui})=\overline{F(p+ui)},$

\item{(b)} the set $P_0:=\{p \in D \mid F(p)=0 \}$ is a 
divide, which is $C^1$ close to
the divide P, and hence the divides $P$ and $P_0$ are combinatorially
equivalent,

\item{(c)} the function $F$ has only non degenerate 
singularities, which are all real,

\item{(d)} for all 
$t \in \bC,\ |t| \leq \eta$ the intersection 
$K_t:=\{(p+iu) \in B_{\lambda} \mid  F(p+iu)=t \} \cap \partial{B_{\lambda}}$ 
is transversal and by a small isotopy
equivalent to the link $L(P)$,

\item{(e)} for the link $K_{\eta}$ the surface 
$\{(p+iu) \in B_{\lambda} \mid F(p+iu)=\eta\}$ 
is  a connected  smooth
symplectic spanning surface of genus $\delta-r+1$ in the $4$-ball.

\endproclaim

\demo {\bf Proof}
Let the divide $P$ be given by smooth
parametrized curves
$\gamma_l:[-1,1] \to \bR^2, 1\leq l \leq r$. Using the 
Weierstrass Approximation Theorem, we can
construct polynomial approximations $\gamma_{l,0}:[-1,1] \to \bR^2$ being 
$C^2$ close to $\gamma_l$ and henceforth give a divide $P_0$ 
with the combinatorics of the divide $P.$ We may choose $\gamma_{l,0}$ such that
$\gamma_{l,0}(s) \notin D, |s| > 1.$ 
Let 
$F:\bC^2 \to \bC$ be a real polynomial map such $F=0$ is a regular
equation for the union of the images of $\gamma_{l,0}.$
Let $S_{\lambda}^3$ be the sphere 
$S_{\lambda}^3:=\{p+iu \in \bC^2 \mid ||p||^2 + \lambda^{-2}||u||^2=1\}.$ For a
sufficiently small $\lambda > 0 $ we have that
the $0$-level of $F$ on $S_{\lambda}^3$ is a model for the link
$L(P).$ For $t \in \bC, t \not=0,$ and $t$ sufficiently small, say $|t| \leq
\eta$, the surface
$X_t:=\{p+iu \in B_{\lambda}^4 \mid F(p+iu)=t \}$ 
is connected and smooth of genus $\delta-r+1$, and 
has a polynomial equation, hence is a symplectic surface in the $4$-ball
$B_{\lambda}$ 
equipped with the standard symplectic structure of $\bC^2.$
The intersection $K_{\eta}:=X_{\eta} \cap \partial{B_{\lambda}}$ is also a 
model for the link $L(P)$  and has hence a symplectic filling with the
required properties.
\enddemo

\remark {\bf Remark} Unfortunately, it is not the case that the restriction of 
$F$ to $B_{\lambda}$ is a 
fibration with only
quadratic singularities, such that for some $\eta > 1$ the fibers
$f_0^{-1}(t), t \in \bC, |t| < \eta,$ are transversal to the boundary of
$B_{\lambda}.$ So, we do not know, as it is the case for divides coming from
plane curve singularities, if it is possible to fill in with a
Picard-Lefschetz fibration, which is compatible with the contact and
symplectic structure.
\endremark

\S 5.  {\bf The gordian number of the link of a divide}\par

The $\delta$-invariant of
a plane curve singularity $S$ is the number of local double 
points in $\bC^2$, that occur in the union of its branches after a small 
generic deformation of the parametrizations 
of the branches. The $\delta$-invariant is also the dimension 
as $\bC$ vector space 
of the quotient of the normalisation of the local ring of $S$ by the local
ring of $S$. The \"Uberschneidungs\-zahl or gordian number $s(L)$
of a link $L$ in $S^3$ is 
the smallest number of cutovers, see  Fig. $5$, by which the link 
can be made trivial [W].  

\midinsert
\cline{\epsffile{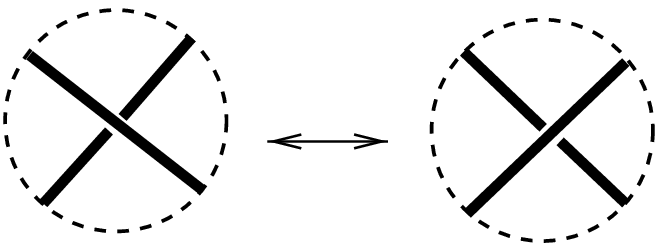}}
\medskip
\centerline{Figure $5$: The move cutover}
\endinsert

J. Milnor proposed the term unknotting number  and
conjectured for  plane curve 
singularities, that the  
unknotting number of the link of the singularity  
equals the $\delta$-invariant of 
the singularity [M]. This conjecture has been proved by P. Kron\-heimer 
and T. Mrowka (see [K1,K2,K3,K-M1,K-M2,K-M3,K-M4,K-M5]). Local 
links of plane curve singularities 
are special among
links of divides and the $\delta$-invariant of a plane curve singularity, which
has only real branches, equals
the number of double points of a divide for the singularity. 
The following theorem extends the computation of the
gordian number to links of divides in general. The $4$-ball genus  
of a link
$L$ in $S^3$ is the minimal genus of a  smooth embedded 
oriented surface in $B^4$
bounded by $L$.

\proclaim{\bf Theorem 4}
Let $P$ be a divide with one branche. 
The gordian number and the $4$-ball genus of the knot $L(P)$
equal the number of double points of the divide $P.$
\endproclaim

The proof will be given at the end of this section.
For the proof of Theorem $4$, we can work in the ball $B_{\lambda}$ 
and use the arguments of Kronheimer and
Mrowka as in their work on the Thom conjecture, together 
with their extension to the  relative
case [K-M3] of a theorem of Taubes [Ta]. In the proof below we will 
use  the global 
curve given by the polynomial equation $\{F=t\}$ of Theorem $3$
together with its completion in $P^2(\bC)$ and apply the affirmative 
answer of  Kronheimer and
Mrowka  to the Thom conjecture [M-K2]. It is also possible to
compute the Thurston-Bennequin number directly from the combinatorial data of
the divide using a global deplacement in the direction of the left 
invariant vectorfield given by $j$ on $S^3$
and to conclude with an inequality of D. Bennequin [E] (see [G]), 
that the
number of crossing points of the divide is a 
lower bound for the gordian number of its link. The point is
that here luckily, in view of the inequality of Bennequin and the lemma below,
the linking number of $L(P)$ with  $j.L(P)$ is maximal among the 
displacements $X.L(P)$ 
given by 
global non-vanishing 
vector fields $X$, which are tangent to the left invariant contact 
distribution spanned by 
$[k,j]$, and therefore yields
the Thurston-Bennequin number. 

\demo {\bf Proof of Theorem $4$}
Let $P$ have $\delta$ double points.
Let $X \subset P^2(\bC)$ be the projective curve given by the equation
$\{F=0\}$ of Theorem $3$.  The curve $X$ 
intersects $\partial{B_{\lambda}}$ transversally 
and has $\delta$ transversal double points in $B_{\lambda}$. Let $Y:=\{F=t\}$
be a non-singular approximation of $X$. Since the genus 
of $Y$ is minimal among all smooth surfaces in $P^2(\bC)$ 
representing $[Y]$ by the work of 
Kronheimer and Mrowka on the Thom conjecture,  the
genus of $Y \cap B_{\lambda}$ is minimal among 
all smooth spanning surfaces in
$B_{\lambda}$ of the link $Y \cap \partial{B_{\lambda}}$ and the $4$-ball
genus of the link $Y \cap \partial{B_{\lambda}}$ equals $\delta$. 
Since the links 
$Y \cap \partial{B_{\lambda}}$ and $L(P)$ are equivalent, we conclude that the
$4$-ball genus of the link $L(P)$ 
equals  $\delta.$ At this point it follows that the 
gordian number of the link $L(P)$
equals or exceeds $\delta$ since the $4$-ball genus of a link is a lower bound
of its gordian number. It is  consequence of the following lemma 
that the gordian number of $L(P)$ equals $\delta$.
\enddemo

\proclaim{\bf Lemma} 
Let $P$ be a divide with $\delta$ double points. The
gordian number of the link $L(P)$ does not exceed $\delta.$
\endproclaim

\demo{\bf Proof}
We will produce an isotopy from the link $L(P)$ to the trivial link 
which has exactly $\delta$ cutovers. First we need to choose a
co-orientation of the branches of the divide. Next we at each point
$p$ of
$P$ we consider 
the normal vector $n_p$ in the direction of the choosen co-orientation such
that for its length we have 
the rule $||p||^2+||n_p||^2=1.$ For $\sigma \in
[0,\pi/2)$ we define $L(P,\sigma)$ to be the link, possibly with
transversal self  crossings,
$$L(P,\sigma):=\{(x,\cos(\sigma)u+\sin(\sigma)n_x)  \in T(D) \mid  (x,u)\in L(P) \} \subset S(T(\bR^2))=S^3$$
The link $L(P,\sigma)$ will have a singularity above the 
crossing point $c$ of the
divide $P$ if $\sigma={\pi-\alpha_c \over 2}$ where $\alpha_c$ is the angle in
between the two normals to $P$ at $c$. The link $L(P,\sigma_{0})$ is trivial if 
$\pi/2 > \sigma_{0} >  {\pi-\alpha_c \over 2}$ for all 
crossing points of $P$, since $L(P,\sigma_{0})$
is spanned by the union of 
embedded disks $\cup_{t \in [\sigma_{0},\pi/2]}L(P,t)\subset S^3$. 
Indeed, observe
that the above formula defines a curve $L(P,\pi/2)$ which  is 
a disjoint union of embedded arcs in $S^3$ and 
that $\cup_{t \in [\sigma_{0},\pi/2]}L(P,t)$ is 
a disjoint union of smoothly embedded $2$-disks in $S^3$. 
The family $L(P,t), t \in [0,\sigma_{0}],$ connects the link $L(P)$ with the
trivial link and has $\delta$ cutovers.

\remark{\bf Example} The knot of the divide hart (Fig. $6$) 
with 2 double points 
is the knot with 10 crossings $10_{145}$ (see
the table F.1 of [Ka] page 261). From Theorem 4 it follows
that the 4-ball genus and the gordian number of the knot $10_{145}$ are 2,
which allows us to complete  entries of  table F.3 of [Ka]. As I learned from 
T. Tanaka he has determined by an 
other method the gordian number of the knot $10_{145}$ [T]. 
The gordian numbers of the
knots $10_{139}$ and $10_{152}$ are proved to be $4$ by Tomomi Kawamura [Kaw]
and she deduced from $s(10_{139})=4$ the gordian number $s(10_{161})=3$.
\endremark

\midinsert
\cline{\epsffile{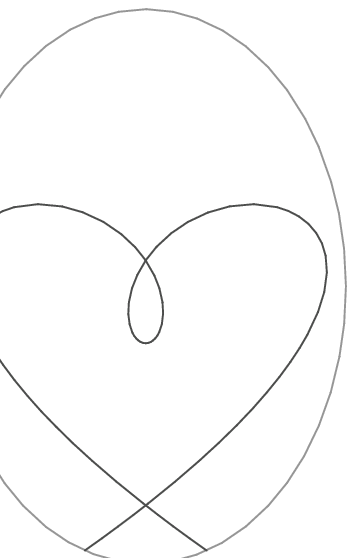} \ \epsffile{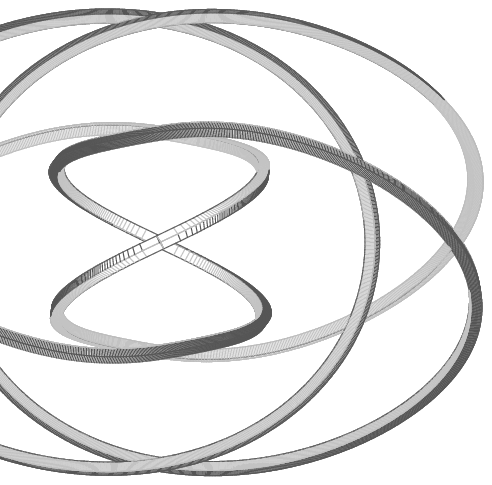}}
\medskip 
\centerline{Figure $6$: The divide hart and its knot $10_{145}$}
\endinsert

\remark{\bf Remark} Let $P$ be a divide with two branches  $P_1$ and $P_2$.
The homological linking number of the two oriented knots $L(P_1)$ and 
$L(P_2)$ equals the number of intersection points of $P_1$ and $P_2$. The
minimal number of cutovers needed to separate by a smooth $2$-sphere the
components of the link $L(P)$ equals also the number of intersection points of
$P_1$ and $P_2$. It follows with Theorem $4$ that the gordian number 
of the link of a divide
equals the number of double points of the divide.
\endremark

\remark{\bf Remark} H. Pinkham proved that for links of singularities of plane
curves the gordian number does not exceed the  
$\delta$-invariant of the singularity. Since the $\delta$-invariant of a
singularity is also the number of double points of a divide of a 
singularity of the same topological 
type, we have reproved the result of Pinkham.
\endremark
 
\par

\bye